\documentclass[12pt]{article}
\usepackage{amssymb,amsmath, amsthm}
\usepackage{graphicx}

\newtheorem{thm}{Theorem}
\newtheorem*{theo}{Theorem}
\newtheorem*{propo}{Proposition}

\newenvironment{defin}{\medskip\noindent{\sc
Definition}.}{\goodbreak\medskip}
\newenvironment{nota}{\medskip\noindent{\sc
Notations}.}{\goodbreak\medskip}
\newenvironment{remk}{\noindent{\sc
Remark}.}{\goodbreak\vskip10pt}

\def\demo{\medskip\goodbreak\noindent
     \hbox{\sc Proof \kern .3em}\ignorespaces}%
  \def \qedbox{$\square$}%
  \def \qed{\hglue1mm\hfill{\ifmmode\qedbox
     \else\unskip\ \hglue0mm\hfill\qedbox\medskip
      \goodbreak\fi}}%

\def\qed{\hglue1mm\hfill\raise -2pt\hbox{\vrule\vbox to 10pt{\hrule width
4pt
                  \vfill\hrule}\vrule}}

\newcommand{\T}{\mathbb {T}}
\newcommand{\A}{\mathbb {A}}

\newcommand{\R}{\mathbb {R}}
\newcommand{\Q}{\mathbb {Q}}

\newcommand{\Z}{\mathbb {Z}}

\newcommand{\Uc}{\mathcal {U}}

\newcommand{\Cc}{\mathcal {C}}

\newcommand{\Gc}{\mathcal {G}}

\newcommand{\Sc}{\mathcal {S}}
\newcommand{\Tc}{\mathcal {T}}

\begin{document}
\title{
Boundaries of instability zones for      symplectic twist maps }
\author{M.-C. ARNAUD
\thanks{ANR DynNonHyp ANR BLAN08-2-313375}
\thanks{Universit\'e d'Avignon et des Pays de Vaucluse, Laboratoire d'Analyse non lin\' eaire et G\' eom\' etrie (EA 2151),  F-84 018Avignon,
France. e-mail: Marie-Claude.Arnaud@univ-avignon.fr}
}
\maketitle
\abstract{   
We construct a $C^2$ symplectic twist map $f$ of the annulus that has an essential invariant curve $\Gamma$ such that:
\begin{enumerate}
\item[$\bullet$] $\Gamma$ is not differentiable;
\item[$\bullet$] the dynamics of $f_{|\Gamma}$ is conjugated to the one of a Denjoy counter-example;
\item[$\bullet$] $\Gamma$ is at the boundary of an instability zone for $f$.

\end{enumerate}
}
 \newpage

\tableofcontents
   \section{Introduction}
  
  The  exact symplectic twist maps of the two-dimensional annulus\footnote{all the definitions are given in subsection \ref{ss12}}  were studied for a long time because they represent (via a
symplectic change of coordinates) the dynamics of the generic symplectic diffeomorphisms of
surfaces near their elliptic periodic points (see \cite{Ch1}). One motivating  example of such  a map was introduced by Poincar\'e for the study of  the restricted 3-Body problem.\\

The   study of such maps was initiated by G.D.~Birkhoff in the 1920s (see \cite{Bir1}). Among other beautiful results, he proved that any essential invariant curve by a symplectic twist map of the annulus is the graph of a Lipschitz map (an {\em essential} curve is a simple loop that is not homotopic to a point). He then introduced the notion of {\em instability zone}.

\begin{defin} An  {\em instability zone} of a symplectic twist map $f$ of the annulus is an open subset $U$ of the annulus $\A$ that is invariant by $f$ and such that:
\begin{enumerate}
\item[--] $U$ is homeomorphic to (open) the annulus $\A$;
\item[--] the closure $\bar U$ of $U$ in $\A$ contains no essential invariant curve that is not contained in the boundary $\partial U$;
\item[--] $U$ is a maximal (for the inclusion $\subset$) subset of $\A$ that satisfies all these properties.
\end{enumerate}
\end{defin}
There are three kinds of instability zones $U$:
\begin{enumerate}
\item  the whole annulus $\A$ can be an instability zone; this happens for example for the standard map with a large enough parameter (see \cite{Au, Mat3, MacPer}); 
\item $U$ is a connected component of the complement of an essential invariant curve;
\item $U$ is bounded; in this case, there exists two Lipschitz  functions $\psi_-<\psi_+$ whose graphs are invariant and that satisfy: $U=\{ (\theta, r)\in\A; \psi_-(\theta)<r<\psi_+(\theta)\}$.
\end{enumerate}
A lot of things are known about the existence of those instability zones.  G.~D.~Birkhoff proved in \cite{Bir2} the existence of such instability zones. He even gave the first (and only) explicit example of boundary for an instability zone. To visualize his example,   imagine the time one map $T$ of the rigid pendulum. It is a symplectic twist  map with one hyperbolic fixed point and two separatrices connecting this fixed point to itself. Perturb $T$ to create one transverse homoclinic intersection at  a point of the upper separatrix without changing the lower separatrix $\Sc$. Then $\Sc$ becomes the boundary of an instability zone that is above $\Sc$ (`above' is `inside' in the following picture).\\
\vspace{.2in}
\begin{center}
\includegraphics[width=9cm]{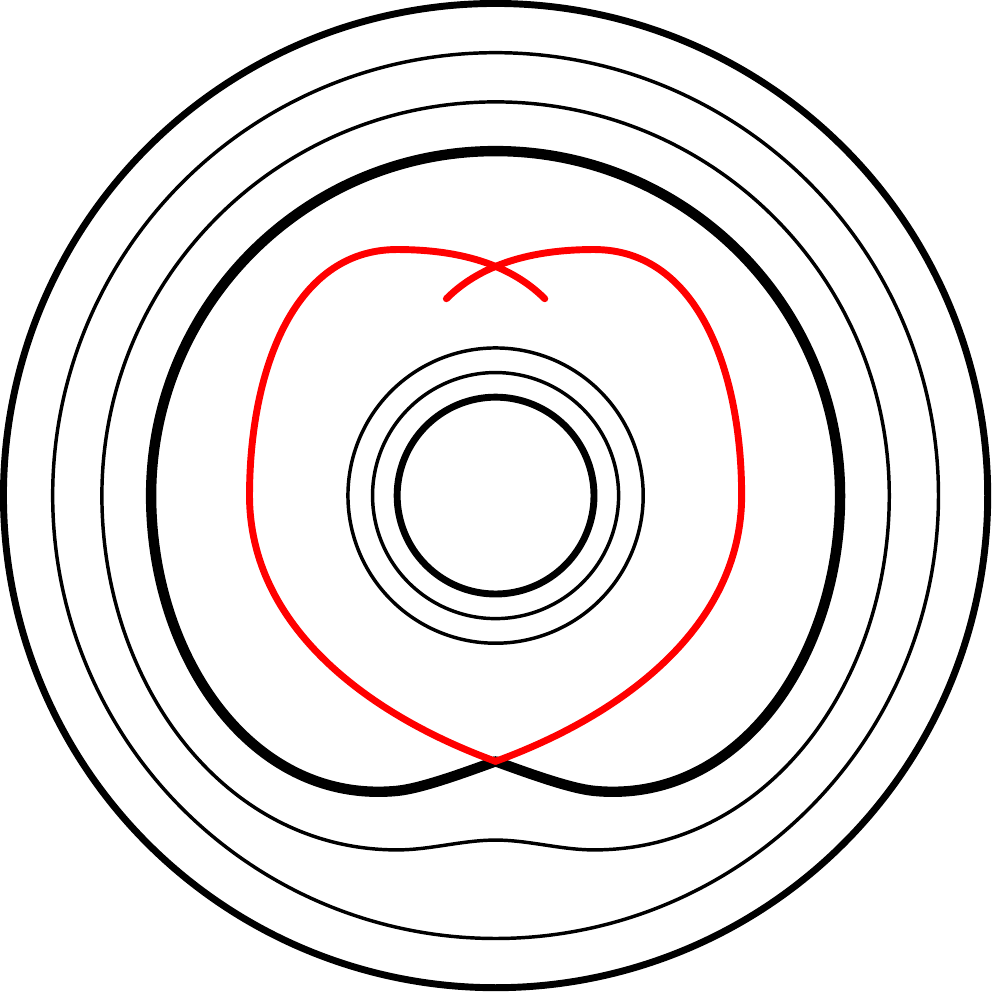}
\end{center}

In Birkhoff's example, the boundary of the instability zone in non-smooth. Modifying the potential in such a way that it has a degenerate minimum, then we obtain an similar example for which the boundary of the instability zone is smooth:\\
\vspace{.2in}
\begin{center}
\includegraphics[width=9cm]{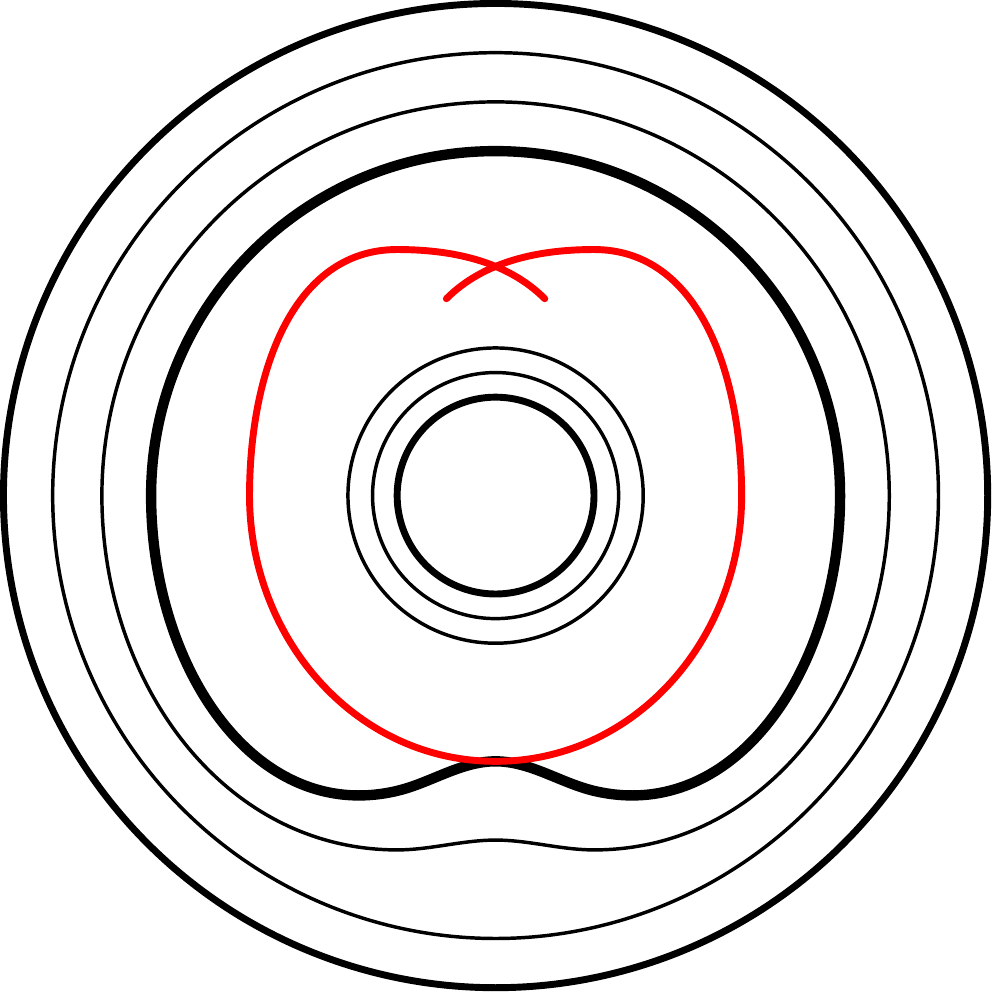}
\end{center}

In \cite{He2}, M.~Herman proved that in general, the boundaries of the instability zones have an irrational rotation number. Hence  Birkhoff's example is not `general'. \\
Curiously, no other examples of explicit boundaries of instability zone are known. To be  complete, let us just mention that in \cite{Mat1}, J.~N.~Mather  proves that the billiard map of a convex billiard whose curvature vanishes at at least one point has an instability zone bounded by the boundary of the billiard phase space. Unfortunately, in this case (vanishing curvature), the billiard map is not a twist map\dots\\

Though we don't know how the `general' boundaries of the instability zones are, we know some facts about what cannot be such a boundary for a sufficiently regular symplectic twist map, for example $C^\infty$:
\begin{enumerate}
\item it cannot be a curve on which the dynamics is $C^\infty$ conjugated to a Diophantine rotation; indeed, KAM theorems (see \cite{Ko, Arno, Mo, Ru, He2} for example) implies that such a curve is accumulated from below and above by other invariant curves;
\item  it cannot be a curve on which the dynamics  is $C^\infty$ conjugated to a rational rotation; indeed, it is proved in the thesis of R.~Douady that in this case you can again apply KAM theory.
\end{enumerate}
Hence a curve that is at the boundary of an instability zone either  is not very regular or has a rational or Liouville rotation number. We then raise the question:\\
{\bf Question.} {\sl Can the boundary of an instability zone with an irrational rotation number be non-differentiable? Can it be smooth?}\\
We will give some answers to these questions in the case of low regularity ($C^1$ or $C^2$). At first, we will prove: 
  \begin{thm}\label{T1}
Let $\omega\in\R\backslash \Q$ be an irrational number. In any neighborhood of $(\theta, r)\rightarrow (\theta+r, r)$ in the $C^2$ topology, there exists a symplectic $C^2$ twist map $f$ of the annulus that has an essential invariant curve $\Gamma$ such that:
\begin{enumerate}
\item[$\bullet$] $f_{|\Gamma}$ is  is $C^0$-conjugated to a Denjoy counter-example and its rotation number is $\omega$;
\item[$\bullet$] if $\gamma~: \T\rightarrow \R$ is the map whose graph is  $\Gamma$, then $\gamma$ is $C^1$ at every point except along the projection of one wandering orbit $(x_n)$, along which $\gamma$ has distinct right and left derivatives;
\item [$\bullet$] $\Gamma$ is the upper boundary of an instability zone $\Uc$ of $f$  ;
\item [$\bullet$] there exists two families of $C^2$ curves  $\gamma_n^s, \gamma_n^u: \R\rightarrow \A$ such that $\gamma_n^u(0)=\gamma_n^s(0)=x_n$ and:
\begin{enumerate}
\item $f\circ \gamma_n^u=\gamma_{n+1}^u$ and $f\circ \gamma_n^s=\gamma_{n+1}^s$;
\item $\forall y\in \gamma_0^s (\R) , \displaystyle{\lim_{n\rightarrow +\infty}d(f^ny, f^nx_0)=0}$ and 
$\forall y\in  \gamma_0^u (\R), \displaystyle{\lim_{n\rightarrow +\infty}d(f^{-n}y, f^{-n}x_0)=0}$;
\item $\gamma_n^s(]-\infty, 0])\cup\gamma_n^u([0, +\infty[)\subset \Gamma$ and $\gamma_n^s(]0, +\infty[)\cup \gamma_n^u(]-\infty, 0[)\subset \Uc$.
\end{enumerate}
\end{enumerate}
\end{thm}
\vspace{.2in}
\begin{center}
\includegraphics[width=11cm]{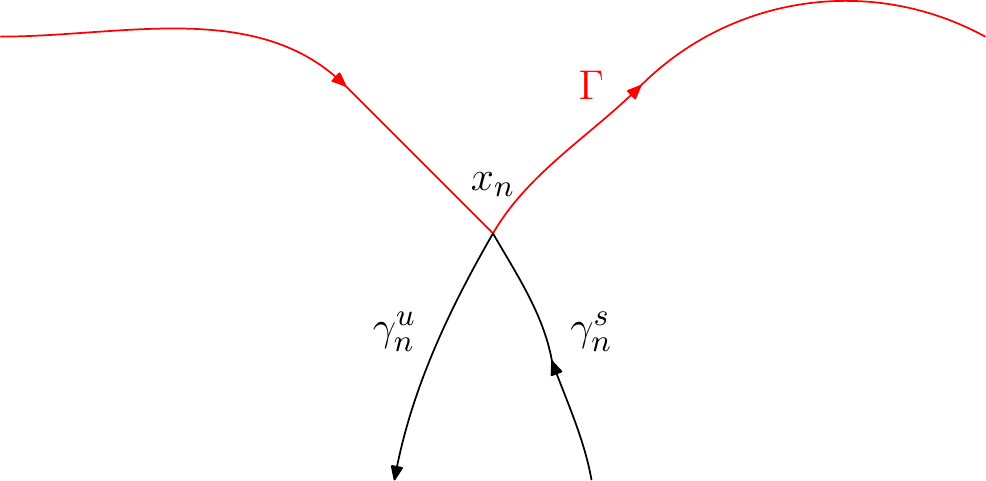}
\end{center}

\begin{remk} 1) With a slight change in the construction, we can ask that $\Gamma$ is the lower boundary of the instability zone $\Uc$;\\
2) if we use a Denjoy counter-example with two disjoint orbits of wandering intervals, we can do the same construction along two orbits $(x_n)$ and $(y_n)$ and the obtain that $\Gamma$ is the common boundary of {\em two} instability zones, the one that is above $\Gamma$ and the other that is under $\Gamma$;\\
3) our counter-example is defined by: $f_\varphi(\theta, r)=(\theta+r, r+\varphi(\theta+r))$ with $\int_\T\varphi=0$. Hence $f_\varphi(\theta, r+1)= f_\varphi(\theta, r)+(0,1)$. If a graph $\Gamma$  is invariant by $f_\varphi$, all the   translated graphs $\Gamma+(0, k)$ with $k\in\Z$ are invariant by $f_\varphi$. This implies that the instability zones of $f_\varphi$ are either the whole annulus $\A$ or bounded instability zones. Hence $\Uc$ is a bounded instability zone, but the theorem gives us the description of just one connected component of its boundary.\\
4) In \cite{Arna2}, we gave an example of a $C^1$ symplectic twist map that has a non differentiable essential invariant curve with irrational rotation number; here, we improve the construction in the following way:
\begin{enumerate}
\item[$\bullet$] using the construction of M.~Herman that is given in \cite{He2}, we manage to improve the regularity of our example and obtain a $C^2$ counter-example;
\item[$\bullet$] using a function $\varphi:\T\rightarrow \R$ whose restriction to a lot of intervals is linear, we manage to create a non-trivial stable manifold for the invariant curve; this and the fact that the rotation number of $f_\varphi$ restricted to the curve is irrational  imply that the invariant curve is at the boundary of an instability zone (see subsection \ref{ss13} for details).\\
Indeed, with the notations of this theorem, $\gamma_n^s(\R)$ is a part of of the stable manifold  of $\Gamma$
$$W^s(\Gamma)=\{ x\in\A; \lim_{n\rightarrow +\infty}d(f^nx, \Gamma)=0\}$$
and 
 $\gamma_n^u(\R)$ is a part of of the unstable manifold  of $\Gamma$
$$W^u(\Gamma)=\{ x\in\A; \lim_{n\rightarrow +\infty}d(f^{-n}x, \Gamma)=0\}.$$
\end{enumerate}
\end{remk}
We now explain how the restricted dynamics to any invariant curve of a symplectic twist map that has an irrational rotation number can become the dynamics at the boundary of an instability zone:

\begin{thm}\label{T2}
Let  $\Gamma$ be an essential invariant curve of a $C^1$ symplectic twist map  $f:\A\rightarrow \A$ whose rotation number is irrational or whose rotation number is rational and the dynamics restricted to $\Gamma$  is:
\begin{enumerate}
\item[--] either $C^0$ conjugated to a rational rotation;
\item[--] or has only hyperbolic periodic points.
\end{enumerate} Then in any neighborhood $\Uc$ of $f$ for the $C^1$ topology, there exists a $C^1$ symplectic twist map $g:\A\rightarrow \A$ such that:
\begin{enumerate}
\item $\Gamma$ is at the boundary of an instability zone of $g$;
\item $g_{|\Gamma}=f_{|\Gamma}$.
\end{enumerate}
\end{thm}
We have seen before that such a result is not valid in $C^\infty$ topology because of the KAM theorems. The tools used to prove   theorem \ref{T2} are specific to the $C^1$ topology: they are the connecting lemma of Hayashi (see \cite{Hay1}) and more precisely some consequences of this connecting lemma that are given in \cite{ABC}.

Contrary to Birkhoff's counter-example or  to theorem \ref{T1}, we have no idea of how the stable/unstable manifold of $\Gamma$ is  in theorem \ref{T2}. Observe too that   for our example of theorem \ref{T1}, we only know a part of the stable/unstable manifold. Hence we raise the following question. \\

\noindent{\bf Question.} {\em  Is it possible to describe (in general or for some specific examples) the stable/unstable manifold of the boundary of an instability zone?}\\
In \cite{LeC1}, P.~Le~Calvez proves interesting facts concerning the topological structure of those sets.\\
We observe in subsection \ref{ss13}  that the existence of a non-trivial stable manifold for an essential invariant curve with an irrational rotation number  implies that this curve is at the boundary of an instability zone. Hence a related question is:\\

\noindent{\bf Question.} {\em Can an irrational essential invariant curve carry a non-uniformly hyperbolic invariant measure?}\\
Indeed, if this happens, the union of the stable and unstable manifold of the invariant measure cannot be contained  in the curve and the curve is then at the boundary of an instability zone.\\

Finally, concerning the first questions that we raised, we obtain an answer just in the case of low regularity. Hence the following questions remain open.\\

\noindent{\bf Questions} {\em 1) Does there exist a smooth ($C^1$, $C^2$, \dots)  curve with an irrational rotation number that is at the boundary of an instability zone for a $C^k$ symplectic twist map with $k\geq 2$?\\
2) Does there exist a non $C^1$  curve with an irrational rotation number that is at the boundary of an instability zone for a $C^k$ symplectic twist map with $k\geq 3$?\\ 
3) How is a `typical' boundary of instability zone (is it regular, how is its rotation number\dots)?
}
\subsection{Structure of the article}

In section \ref{S2}, we will recall Michel Herman's construction of a $C^2$ symplectic twist map that has an essential invariant curve on which the dynamics is Denjoy. In particular, we will give some useful estimates in subsection \ref{ss23}.

In section \ref{S3} we will build the counter-example that is descibed in theorem \ref{T1}. In subsection \ref{ss31}, we will construct the homeomorphism that will represent the projected dynamics along the invariant curve and we will be more precise about the choice of the constants in subsection \ref{ss32}. In subsections \ref{ss33} and \ref{ss34}, we will prove some estimates. Then, in subsection \ref{ss35}, we will prove that our modified example is $C^2$ and in section \ref{ss36} we will determine part of the stable/unstable manifolds of the invariant curve.

Finally, we will prove theorem \ref{T2}  in section \ref{S4}.
\subsection{Notations and definitions}\label{ss12}

\begin{nota}
\noindent $\bullet$ $\T=\R/\Z$ is the circle.

\noindent $\bullet$ $\A=\T\times \R$ is the annulus and an element of $\A$ is denoted
by $(\theta, r)$.

\noindent $\bullet$ $\A$ is endowed with its usual symplectic form, $\omega=d\theta\wedge dr$ and its usual Riemannian metric.

\noindent  $\bullet$ $\pi: \T \times \R \rightarrow\T$ is the first projection and $\tilde\pi:
\R^2\rightarrow \R$ its lift.

\end{nota}

\begin{defin} A $C^1$ diffeomorphism $f: \A\rightarrow \A$ of the annulus that is isotopic
to identity  is a {\em positive twist map} (resp. {\em negative twist map}) if, for any given lift $\tilde f: \R^2\rightarrow
\R^2$ and for every
$\tilde\theta\in\R$, the maps $r\mapsto \tilde\pi\circ \tilde f(\tilde\theta,r)$ is an increasing (resp decreasing)   diffeomorphisms. A {\em
twist map} may be positive or negative. 
\end{defin}
Then the maps $f_\varphi$ that we defined just after theorem \ref{T1}  are positive symplectic twist maps.

\begin{defin} Let $\gamma~:\T\rightarrow\R$ be a continuous map. We say that $\gamma$ is $C^1$ at $\theta\in \T$ if there exists a number $\gamma'(\theta)\in\R$ such that, for every sequences $(\theta^1_n)$ and $(\theta_n^2)$ of points of $\T$ that converge to $\theta$ such that $\theta_n^1\not=\theta_n^2$, then:
$$\lim_{n\rightarrow \infty} \frac{\gamma(\theta_n^1)-\gamma(\theta_n^2)}{\theta^1_n-\theta^2_n}=\gamma'(\theta) $$
where we denote by $\theta_n^1-\theta_n^2$ the unique number that represents $\theta_n^1-\theta_n^2$ and that belongs to $]-\frac{1}{2}, \frac{1}{2}]$.
\end{defin}
If we assume that $\gamma$ is differentiable at every point of $\T$, then this notion of $C^1$ coincides with the usual one (the derivative is continuous at the considered point).
\subsection{Stable manifold of   invariant curves} \label{ss13}

A consequence of a theorem of J.~N.~Mather  is that  if an essential curve $\Gamma$ that is invariant by a symplectic twist map is at the boundary of an instability zone, then $W^s(\Gamma)\backslash \Gamma\not=\emptyset$ and this is equivalent too to $W^u(\Gamma)\backslash \Gamma\not=\emptyset$. \\
More precisely, in \cite{Bir1}, G.~D.~Birkhoff proved that if $U$ is an instability zone, if $U_1$ is a neighborhood of one of its ends (i.e, eventually after compactification,  a connected component of its boundary) and $U_2$ is a neighborhood of the the other end, then there exists an orbit  traveling from $U_1$ to $U_2$. This theorem was improved in \cite{Mat2} by J.~Mather who proved that if $\Cc_1$, $\Cc_2$ are the ends of $U$, there exists an orbit whose $\alpha$-limit set is in $\Cc_1$ and $\omega$-limit set is in $\Cc_2$. J.~N.~Mather used variational arguments and after that, P.~Le~Calvez gave in \cite{LeC1} a purely topological proof of this result.

Conversely, let us assume that $\Gamma$ is   an essential invariant curve that is invariant by a symplectic twist $f:\A\rightarrow \A$ and that $W^s(\Gamma)\backslash\Gamma\not=\emptyset$. The example of the rigid pendulum proves that it can happen that $\Gamma$ is not at the boundary at an instability zone\footnote{let us recall that we asked that an instability zone is homeomorphic to the open annulus}. Let us assume that $f_{|\Gamma}$ has an irrational rotation number or that $f_{|\Gamma}$ is $C^0$ conjugated to a rational rotation. Suppose that $\Gamma$ is not at the boundary of an instability zone. Then there exists two sequences of essential invariant curves $(\Gamma_n^-)$ and $(\Gamma_n^+)$ that are different from $\Gamma$ such that:
\begin{enumerate}
\item $\forall n, \Gamma_n^+\geq \Gamma\quad{\rm and}\quad \Gamma_n^-\leq \Gamma$;
\item\label{pt2}  $\displaystyle{ \lim_{n\rightarrow \infty} d(\Gamma_n^-, \Gamma)=0}$ and  $\displaystyle{ \lim_{n\rightarrow \infty} d(\Gamma_n^+, \Gamma)=0}$.
\end{enumerate}
Birkhoff's theorem implies that the curves $\Gamma_n^-$, $\Gamma_n^+$ are equi-Lipschitz and then relatively compact for the $C^0$ norm (we speak of the $C^0$ norm of the function  whose graph is the   curve of interest). Let $\Gamma^*$ be any limit point of one of these two sequences. Then $\Gamma^*$ in an essential invariant curve such that $\Gamma\cap\Gamma^*\not=\emptyset$. Hence $f_{|\Gamma^*}$ has the same rotation number as $f_{|\Gamma}$. \\
M.~Herman proved in \cite{He2} that  two curves with the same irrational rotation number are equal. Moreover, if the restriction of a symplectic twist map $f$ restricted to an essential invariant curve $\Gamma$ is $C^0$ conjugated to a rational rotation, all the orbits are action minimizing (see e.g. \cite{Gol}) and a consequence of the results of G.~Forni \& J.~Mather contained in     \cite{FoMa} (see theorem 13.3) is that when an essential invariant curve is filled by periodic orbit, there exists no other minimizing orbit with the same rotation number and then no other invariant curve with the same rotation number. Hence any  other invariant curve  and $\Gamma$ are disjoint. \\ 
Hence $\Gamma^*=\Gamma$ and the two sequences 
$(\Gamma_n^-)$ and $(\Gamma_n^+)$ converge to $\Gamma$. If $\Gamma^\pm_n$ is the graph of $\gamma_n^\pm$, this implies that the sets $\{ (\theta, r); \gamma_n^-(\theta)<r<\gamma_n^+(\theta)\}$ are  a base of neighborhood of $\Gamma$. Because they are invariant by $f$, this implies that $W^s(\gamma)=\Gamma=W^u(\Gamma)$.  We now summarize this result and Mather's result: 
\begin{propo}
Let $f:\A\rightarrow \A$ be a symplectic twist map and let $\Gamma$ be an essential invariant curve. Then:
\begin{enumerate}
\item if $\Gamma$ is at the boundary of an instability zone, then $W^s(\Gamma)\backslash\Gamma\not=\emptyset$ and $W^u(\Gamma)\backslash\Gamma\not=\emptyset$;
\item if $W^s(\Gamma)\backslash\Gamma\not=\emptyset$ or $W^u(\Gamma)\backslash\Gamma\not=\emptyset$ and if the rotation number of $f_{|\Gamma}$ is irrational or if it is rational and if $f_{|\Gamma}$ is $C^0$ conjugated to a rational rotation, then $\Gamma$ is at the boudary of an instability zone.
\end{enumerate}
\end{propo}

  \section{An example due to Michel Herman}\label{S2}
  In \cite{He2}, M.~Herman gives an example of a $C^2$ symplectic twist map $f:\A\rightarrow \A$ that has a $C^1$ invariant curve $\Cc$ such that $F_{|\Cc}$ is $C^0$-conjugated to a Denjoy counter-example. Let us recall his construction. We fix $\omega\in\R\backslash \Q$.
  \subsection{Generalized standard map}
The following  family of symplectic twist maps  was introduced by M.~Herman in \cite{He2}. The  maps are defined  by~:
$$f_\varphi: \T\times\R\rightarrow \T\times\R; (\theta, r)\mapsto (\theta+r, r+\varphi(\theta +r)).$$
where $\varphi~: \T\rightarrow \R$ is a $C^1$ map such that $\int_{\T}\varphi(\theta)d\theta=0$.\\
As noticed by M.~Herman, the main advantage of this map  is the following one.  Using the explicit formula of $f_\varphi$, it is easy to see that the graph of $\psi~: \T\rightarrow \R$ is invariant by $f_\varphi$ if and only if: 
$$\forall \theta\in\T, (\theta +\psi (\theta), \psi (\theta)+\varphi (\theta +\psi(\theta)))=(\theta+\psi (\theta), \psi(\theta+\psi (\theta))).$$
If we rewrite this equality and we denote a lift of $g:\T \rightarrow \T$   by $\tilde g:\R\rightarrow\R$, we obtain the following criterion for the invariance of the graph of $\psi$. The graph of $\psi~: \T\rightarrow \R$ is invariant by $f_\varphi$ if and only if we have:
\begin{enumerate}
\item[$\bullet$] $g=Id_{\T}+\psi$ is an orientation preserving homeomorphism of $\T$;
\item[$\bullet$] $Id_\R +\frac{1}{2} \varphi=\frac{1}{2}\left( \tilde g+ \tilde g ^{-1}\right)$.
\end{enumerate}
In this case, $g$ is none other than  the projected dynamics of the restricted to the graph of $\psi$  dynamics. In particular, the restricted dynamics is conjugated to $g$ (via the first projection).\\

Let us give the idea of the construction of M.~Herman: he builds a particular Denjoy counter-example $g: \T\rightarrow \T$ of rotation number $\omega$. Because of Denjoy's theorem, such a $g$ cannot be $C^2$. By using very clever estimates, M. Herman proves that $\varphi=\tilde g+\tilde g^{-1}-2{\rm Id}$ is $C^2$. Hence $f_\varphi$ is the wanted counter-example.
 
  \subsection{Explicit construction of a circle diffeomorphism }
  We use the construction that is described in \cite{He2} p. 94, with   only a slight change:  we define the function $\eta$ in such a way that the Denjoy counter-example is linear on some small segments. 
  
  Let us recall that we fixed   $\omega\in\R\backslash \Q$. Let us fix  $\delta>0$ and  $C>>1$. \\We introduce for $k\in\Z$:
$$\ell_k=\frac{a_C }{(|k|+C)(\log(|k|+C))^{1+\delta}}$$
where $a_C $ is chosen such that $\displaystyle{ \sum_{k\in\Z}\ell_k=1   }$. \\
We use a $C^\infty$ function  $\eta~:\R\rightarrow \R$ such that $\eta\geq 0$, ${\rm support}(\eta)\subset [\frac{1}{4}, \frac{3}{4}]$, $\eta_{|[\frac{3}{8}, \frac{5}{8}]}=1$, $\eta (1-t)=\eta (t)$  and $\int_0^1\eta(t)dt=1$. We define $\eta_k$ by~: $\eta_k(t)=\eta\left( \frac{t}{\ell_k}\right)$. Then we have: $\int_0^{\ell_k}\eta_k(t)dt=\ell_k$. Moreover, there exist two constants $C_1$, $C_2$, that depend  only on $\eta$, such that~:
\begin{equation} C_1\leq |\eta_k|\leq C_2;\quad \frac{C_1}{\ell_k}\leq |\eta_k'|\leq \frac{C_2}{\ell_k}.\end{equation}We assume now that $C>>1$ is large enough so that:
$$\forall k\in\Z, \left| \frac{\ell_{k+1}}{\ell_k}-1\right|C_2<1.$$
Then the map $h_k~: [0, \ell_k]\rightarrow [0, \ell_{k+1}]$ defined by $h_k(x)=\int_0^x \left(1+\left(\frac{\ell_{k+1}}{\ell_{k}}-1\right) \eta_k(t)\right)dt$ is a $C^\infty$ diffeomorphism such that $h_k(\ell_k)=\ell_{k+1}$.

There exists a Cantor subset $K\subset \T$ that has zero Lebesgue measure and that is such that the connected components of $\T\backslash K$, denoted by $(I_k)_{k\in\Z}$, are on $\T$ in the same order as the sequence $(k\omega)$ and such that ${\rm length}(I_k)=\ell_k$.\\
Let us recall what is the semi-conjugation $j:\T\rightarrow \T$ of the Denjoy counter-example to the  rotation $R_\omega$. If $x\in \{ k\omega; k\in\Z\}$, then we define~: $j^{-1}(x)=\int_0^xd\mu(t)$ where $\mu$ is the probability measure $\displaystyle{\mu =\sum_{k\in \Z}\ell_k\delta_{k\omega}}$, $\delta_{k\omega}$ being the Dirac mass at $k\omega$. Then $j:\T\rightarrow\T$ is a continuous map with degree 1 that preserves the order on $\T$ and that is such that $j(I_k)=k\omega$.\\
Then there is  a $C^1$ diffeomorphism $g:\T\rightarrow \T$ that fixes $K$, is such that $K$ is the unique minimal subset for $g$, has for rotation number $\rho(g)=\omega$, verifies $j\circ g=R_\omega\circ j$. If $k\in\Z$, we introduce the notation: $g_{|I_k}=g_k$; then we have: $g_k(I_k)=I_{k+1}$. Following \cite{He2} again, we can assume that: $g_k'=g'_{|I_k}=\left(1+\left(\frac{\ell_{k+1}}{\ell_{k}}-1\right) \eta_k\right)\circ R_{-\lambda_k}$ where $R_{-\lambda_k}(I_k)=[0, \ell_k]$  and that $g_k: I_k\rightarrow I_{k+1}$ is defined by~: $g_k=R_{\lambda_{k+1}}\circ h_k\circ R_{-\lambda_k}$.\\

  \subsection{Some useful inequalities}\label{ss23}
  We recall without proof some inequalities that are given in \cite{He2} (sometimes we give some slight improvement of these inequalities) and that are useful to prove that $\tilde g+\tilde g^{-1}$ is $C^2$. The constants $C_i$  are independent of $k$ and $C>>1$ and the limits are uniform in $C>>1$.
Introduce the notation: $K_k=\frac{\ell_{k+1}}{\ell_k}-1$. Then  for $C>>1$ large enough:
\begin{equation}\label{E2}{\rm if}\quad n\geq 1, K_{\pm n}=\pm\frac{-1}{n+C}+\frac{\varepsilon(\pm n, C)}{(n+C)^2}\pm \frac{-(1+\delta)}{(n+C)\log(n+C)}
\end{equation}
where $\sup\{ |\varepsilon(\pm n,k)|; n, C\geq 1\}=c<+\infty$. If $n=0$, the good formula is the formula with a +.
 \begin{equation}\label{E3} \frac{C_1}{|k|+C}\leq |K_k|\leq \frac{C_2}{|k|+C} \end{equation}
\begin{equation}\label{E4}C_1K_k^2\leq   |K_{k-1}-K_{k}|=K_{k}-K_{k-1}\leq C_2 K_k^2 \end{equation} 
\begin{equation}\label{E5} \frac{C_2(\log C)^\delta}{(|k|+C)(\log(|k|+C))^{1+\delta}}\geq \ell_k\geq \frac{C_1}{(|k|+C)(\log(|k|+C))^{1+\delta}}
\end{equation}
\begin{equation} \lim_{k\rightarrow \pm\infty}\frac{K_k^2}{\ell_k}=0. \end{equation}
 We don't recall here how we can deduce the fact that $\tilde g+\tilde g^{-1}$ is $C^2$ from these inequalities, because we will give a very similar proof for the modified example in the next section.  
 
 Let us just notice the following fact that is due to our modification of the function $\eta$:
 $$\forall t\in [\frac{3}{8}\ell_k, \frac{5}{8}\ell_k], h_k(t)= \frac{\ell_{k+1}}{\ell_k}t.$$
 
 Let us now give some estimates that were not given in \cite{He2}. We introduce the notation: $m_k:= 1+K_k+\frac{1}{1+K_{k-1}}=\frac{\ell_{k+1}}{\ell_k}+\frac{\ell_{k+1}}{\ell_{k+2}}$. We have: $$m_{k+1}-2-(K_{k+1}-K_k)= \frac{K_k^2}{1+K_k} $$
 hence we deduce from (\ref{E4}) that \begin{equation}\label{E7} |m_{k+1}-2|\leq C_2K_k^2\quad{\rm and}\quad |m_k-m_{k+1}|\leq C_2 K_k^2.\end{equation}  Because of (\ref{E3}) we deduce that: 
    \begin{equation}\label{E8}|m_{k+1}-2|\leq \frac{C_2}{(|k|+C)^2}.\end{equation}

  \section{Modification of Michel Herman's example}\label{S3}
  \subsection{Explicit construction of a circle homeomorphism}\label{ss31}
  We introduce two new functions $\gamma_-, \gamma_+:\R\rightarrow \R$ such that:
  \begin{enumerate}
  \item[$\bullet$] ${\rm support}(\gamma_\pm)\subset [0, 1]$; 
  \item[$\bullet$] $\gamma_{\pm|\R\backslash \{ \frac{1}{2}\}}$ is $C^\infty$;
  \item[$\bullet$] $\gamma_{-|[\frac{1}{2}, 1]}=0$; $\gamma_{+|[0, \frac{1}{2}]}=0$;
  \item [$\bullet$] $\forall t\in [\frac{3}{8}, \frac{1}{2}[, \gamma_-(t)=1$ and 
  $\forall t\in ]\frac{1}{2}, \frac{5}{8}], \gamma_+(t)=1$;
  \item[$\bullet$]  $\int_0^1\gamma_\pm(t)dt=0$.

  \end{enumerate}
  Hence these two functions are $C^\infty$ on $\R\backslash \{ \frac{1}{2}\}$ and discontinuous at the point  $\frac{1}{2}$.\\
  We define a sequence of functions $(\gamma_k)$ by: 
  $$\gamma_k(x)=\gamma_+(\frac{x}{\ell_k})\quad{\rm if}\quad k\geq 1\quad{\rm and}\quad \gamma_k(x)=\gamma_-(\frac{x}{\ell_k})\quad{\rm if}\quad k\leq 0$$

  Let us  fix a sequence $(\alpha_k)$ of real numbers such that $0<| \alpha_k|\leq A.|K_k|$ (where $A$ is a constant). Then we define $\psi_k: \R\rightarrow \R$ by $\psi_k(x)=K_k\eta_k(x)+\alpha_k \gamma_k(x)$ and a new function $h_k: [0, \ell_k]\rightarrow \R$ by $h_k(x)=\int_0^x (1+\psi_k(t))dt$. If $C$ is large enough ($C$ was the constant that is used to define $(\ell_k)$ and then $(K_k)$), then $(K_k)$ and $(\alpha_k)$ are small enough  ($A$ is a fixed constant that doesn't depend on $C$ ) and $1+\psi_k$ is positive. Hence $h_k$ is a homeomorphism onto $[0, \ell_{k+1}]$.\\
  Let us notice that $h_k$ is differentiable everywhere except at $\frac{\ell_k}{2}$ where it has distinct left and right derivatives. More precisely:
 
 1) if $k\geq 1$, then: $\forall x\in [\frac{3}{8}\ell_k, \frac{1}{2}\ell_k], h_k(t)=\frac{\ell_{k+1}}{\ell_k}t$ and 
 $\forall x\in [\frac{1}{2}\ell_k, \frac{5}{8}\ell_k], h_k(t)=(\frac{\ell_{k+1}}{\ell_k}+\alpha_k)t-\frac{\alpha_k\ell_k}{2}$;
 
 2) if $k\leq 0$, then: $\forall x\in [\frac{3}{8}\ell_k, \frac{1}{2}\ell_k], h_k(t)=(\frac{\ell_{k+1}}{\ell_k}+\alpha_k)t-\frac{\alpha_k\ell_k}{2}$ and 
 $\forall x\in [\frac{1}{2}\ell_k, \frac{5}{8}\ell_k], h_k(t)=\frac{\ell_{k+1}}{\ell_k}t$.
 
 Then with this new functions $h_k$ we can construct $g_k$ and $g$ exactly as this was done in Herman's example. The only difference is that there is a discontinuity of $g'$ at the middle of every connected component of the wandering set, the map $g$ being linear on a right neighborhood and on a left neighborhood of each such singularity. \\
 Moreover, $h_k'$ tends to $1$ when $k$ tends to $\pm\infty$. This implies (a precise proof was given in \cite{Arna2}) that  $g$ and the curve $\Gamma$ are $C^1$ at all the points that are not  at the middle of every connected component of the wandering set. Observe that the set of discontinuities of $g$ corresponds to one orbit.
  
 \subsection{Choice of a ``good'' sequence $(\alpha_k)$. }\label{ss32}
 Let us recall that we want that $\varphi=\tilde g+\tilde g^{-1}-2{\rm Id}_\R$ is $C^2$. We need to choose carefully the sequence $(\alpha_k)$ to obtain that. Let us now explain how we choose $(\alpha_k)$, and after that we will prove that $\varphi$ is $C^2$.
 
 We begin by choosing two small $\alpha_1>0$ and $\alpha_0<0$ such that:
 $$\frac{1}{1+K_0+\alpha_0}+1+K_1=\frac{1}{1+K_0}+1+K_1+\alpha_1.$$
 We denote by $m$ this quantity.\\
 Then we extend this sequence by using the constants $m_k:= 1+K_k+\frac{1}{1+K_{k-1}}$:
\begin{equation}\label{E9}\forall k\in\Z\backslash\{0\}, 1+K_{k+1}+\alpha_{k+1}+\frac{1}{1+K_k+\alpha_k}=m_{k+1}.\end{equation}
 If we denote by $\Phi_k$ the map $\Phi_k: ]0, +\infty[\rightarrow \R$ defined by $\Phi_k(t)=m_k-\frac{1}{t}$, each $\Phi_k$ is increasing and we have: $\Phi_{k+1}(1+K_k)=1+K_{k+1}$. Because $\alpha_1>0$, we deduce that we can define $(\alpha_n)_{n\geq 1}$ by using (\ref{E9}) and that: $\forall n\geq 1, \alpha_n>0$. In a similar way, each $\Phi_k^{-1}$ is increasing on  $]-\infty, m_k[$ and $\alpha_0<0$, hence we can define $(\alpha_{-n})_{n\geq 1}$ by (\ref{E9}) and we have then: $\forall n\geq 0, \alpha_{-n}<0$. Similar remarks were done in \cite{Arna2}. \\
 For this particular choice of $(\alpha_k)$, we can notice that for all $k\in\Z$,  $h_k+h_{k-1}^{-1}$ is linear in the interval $[\frac{3}{8}\ell_k, \frac{5}{8}\ell_k]$. More precisely (we use the fact that the $h_k$ are continuous at $\frac{\ell_k}{2}$ to determine some constants)):
 \begin{enumerate}
 \item[$\bullet$] if $k\geq 2$: if $x\in [\frac{3}{8}\ell_k, \frac{1}{2}\ell_k]$, $h_k(x)+h_{k-1}^{-1}(x)=(1+K_k)x+\frac{1}{1+K_{k-1}}x=m_kx$ and if $x\in [\frac{1}{2}\ell_k, \frac{5}{8}\ell_k]$:\\
  $h_k(x)+h_{k-1}^{-1}(x)=(1+K_k+\alpha_k)x-\frac{\alpha_k\ell_k}{2}+\frac{1}{1+K_{k-1}+\alpha_{k-1}}(x+\frac{\alpha_{k-1}\ell_{k-1}}{2})=m_kx$;
 \item[$\bullet$] if $k=1$:  if $x\in [\frac{1}{2}\ell_k, \frac{5}{8}\ell_k]$: \\
  $h_1(x)+h_{0}^{-1}(x)=(1+K_1+\alpha_1)x-\frac{\alpha_1\ell_1}{2}+\frac{1}{1+K_{0}}x=m x-\frac{\alpha_1\ell_1}{2}=m_1x-\frac{\alpha_1\ell_1}{2}$ and  if $x\in [\frac{3}{8}\ell_1, \frac{1}{2}\ell_1]$:\\
   $h_1(x)+h_{0}^{-1}(x)=(1+K_1)x+\frac{1}{1+K_0+\alpha_0}(x +\frac{\alpha_0\ell_0}{2})= m_1x-\frac{\alpha_1\ell_1}{2}$   (let us notice that we change the notation for $m_1$ from this point);
 \item[$\bullet$]   if $k\leq 0$: if $x\in [\frac{1}{2}\ell_k, \frac{5}{8}\ell_k]$, $h_k(x)+h_{k-1}^{-1}(x)=(1+K_k)x+\frac{1}{1+K_{k-1}}x=m_kx$ and    if $x\in [\frac{3}{8}\ell_k, \frac{1}{2}\ell_k]$:\\
  $h_k(x)+h_{k-1}^{-1}(x)=(1+K_k+\alpha_k)x-\frac{\alpha_k\ell_k}{2}+\frac{1}{1+K_{k-1}+\alpha_{k-1}}(x+\frac{\alpha_{k-1} \ell_{k-1}}{2})=m_kx$.
 \end{enumerate}
 We deduce immediately that the function $\varphi=\tilde g+\tilde g^{-1}-2{\rm Id}_\R$ is linear on each segment $J_k\subset I_k$ that is at the middle of $I_k$ and has length $\frac{\ell_k}{4}$. In particular, the restriction of $\varphi$ to the interior of any interval $I_k$ is $C^\infty$.
 
 We   denote by $\varphi_k$ the $C^\infty$ function that is equal to $\varphi$ on $I_k$ and equal to $0$ everywhere else. Then: $\varphi=\sum\varphi_k$ and to prove that $\varphi$ is $C^2$, we just have to prove that $\displaystyle{\lim_{k\rightarrow\pm\infty} \|ÊD^2\varphi_k\|_{C^0}=0}$. If we want to prove that $\varphi_k$ is  close to 0 in $C^2$ topology, we have to prove that $\displaystyle{\lim_{C\rightarrow+\infty}\sup\{  \|ÊD^2\varphi_k\|_{C^0}; k\in\Z\}=0}$.
  
  \subsection{Estimation of $(\alpha_n)_{n\geq 1}$}\label{ss33}
  If we want to have a control on $\|ÊD^2\varphi_k\|_{C^0}$, we need to have a control of the sequence $(\alpha_k)$.  We use the following notation: $\beta_k=K_k+\alpha_k$.\\
  We have built the sequences $(\ell_k)$, $(K_k)$ and $(m_k)$ that depend on a certain constant $C>>1$, we have chosen $\alpha_1>0$ small and defined: 
  $$\forall n\geq 1, 1+\beta_{n+1} +\frac{1}{1+\beta_n}=m_{n+1}.$$

We have considered the functions $\Phi_k: ]0, +\infty[\rightarrow \R$ defined by $\Phi_k(t)=m_k-\frac{1}{t}$.  Then we have: $1+  \beta_{k+1}=\Phi_{k+1}(1+\beta_k)$. This function is strictly increasing and concave. When  $m_k>2$, $\Phi_k$ has exactly two fixed points $a_k<1<b_k$ and we have: $b_k=\frac{1}{2}(m_k+\sqrt{m_k^2-4})$ hence (see  (\ref{E8})): \begin{equation}\label{E10}  0<b_k-1<  C_1\sqrt{m_k-2}\leq \frac{C_2}{n+C} \end{equation}

Let us now compare $1+  \beta_{n+1}=\Phi_{n+1}(1+\beta_n)$ with $  1+\beta_n$. We fix  a constant $B>>2$. There are three cases:
\begin{enumerate}
\item if $m_{n+1}\leq 2$ and    $\beta_n\leq \frac{B}{n+1+C}$; then $1+\beta_{n+1}=\Phi_{n+1}(1+\beta_n)\leq \Phi_{n+1}(1+ \frac{B}{n+1+C})\leq 1+ \frac{B}{n+1+C}$ because   $\Phi_{n+1}\leq{\rm Id}$;
\item \label{pt2} if  $m_{n+1}>2$ and $\beta_n\leq \frac{B}{n+1+C}$; then $1+\beta_{n+1}=\Phi_{n+1}(1+\beta_n)\leq \Phi_{n+1}(1+ \frac{B}{n+1+C})\leq 1+ \frac{B}{n+1+C}$ because $\Phi_{n+1 |[b_{n+1}, +\infty[}\leq {\rm Id}_{|[b_{n+1}, +\infty[}$ and $1+ \frac{B}{n+1+C}=1+ \frac{C_2}{n+1+C}+ \frac{B-C_2}{n+1+C}\geq b_{n+1}$ if $B$ is large enough (see (\ref{E10})); 
\item  if  $1+\beta_n> 1+\frac{B}{n+1+C}$. We introduce the notation $\delta_n=\frac{B}{2(n+1+C)}$. The function $\Phi_{n+1}$ being concave such that $D\Phi_{n+1}(1+\delta_{n+1})=\frac{1}{(1+\delta_{n})^2}$, we have:
$$1+\beta_{n+1}-\Phi_{n+1}(1+\delta_n)\leq \frac{1}{(1+\delta_{n})^2}(1+\beta_n-(1+\delta_{n})).$$
If $m_{n+1}\leq 2$, we have  $\Phi_{n+1}(1+\delta_n)\leq1+\delta_n$ because $\Phi_{n+1}\leq {\rm Id}$; if $m_{n+1}>2$, as
 $1+\delta_n>b_{n+1}$ (see point \ref{pt2}), we have $\Phi_{n+1}(1+\delta_n)\leq1+\delta_n$ and then:
$$ \beta_{n+1} \leq (1-\frac{1}{(1+\delta_{n})^2}) \delta_n+\frac{1}{(1+\delta_{n})^2}  \beta_n.$$
As $\delta_n\leq \frac{\beta_{n}}{2}$, we deduce: $$\beta_{n+1}\leq \left(\frac{1}{2} + \frac{1}{2}  \frac{1}{(1+\delta_{n})^2}    \right) \beta_n=\left(\frac{1}{2} + \frac{1}{2}  \frac{1}{(1+\frac{B}{2(n+1+C)})^2}    \right) \beta_n.$$
We deduce for $C$ large enough:
$$\beta_{n+1}\leq (1-\frac{ B}{3(n+C+1)})\beta_n.
$$
\end{enumerate}
We choose $B\geq 3$. We have then:
\begin{equation}\label{E11}\beta_{n+1}\leq (1-\frac{1}{n+C+1})\beta_n=   \frac{n+C }{n+C+1}\beta_n.
\end{equation}

 Let us now prove some estimates for $(\beta_n)$ (and then $(\alpha_n)$). At first, let us recall that $\beta_n>K_n$ (because we have noticed that $\alpha_n>0$). Let us now choose  $\alpha_1>0$ small enough such that $\beta_1=\alpha_1+K_1 \leq\frac{B}{1+C}$; this is possible because $K_1<0$ (see (\ref{E2})). Now we prove by recurrence that:
$ \forall n\geq 1, \beta_n\leq \frac{B}{n+C}$.\\
The result is true for $n=1$. \\
Let us assume that it is true for some $n\geq 1$. There are two cases:
\begin{enumerate}
\item[$\bullet$] either $\beta_n\leq \frac{B}{ n+1+C}$; then we have proved that: $\beta_{n+1}\leq \frac{B}{n+1+C}$;
\item[$\bullet$] or $\beta_n>  \frac{B}{ n+1+C}$; then by (\ref{E11}), we have: $\beta_{n+1}\leq \frac{n+C}{n+1+C}\beta_n\leq \frac{n+C}{n+1+C} \frac{B}{n+C}= \frac{B}{n+1+C}$.
\end{enumerate}
Finally, we have proved that:
$$\forall n\geq 1, K_n\leq \beta_n\leq  \frac{C_2}{n+C}.$$
Using (\ref{E3}), we deduce similar estimates for $\alpha_n=\beta_n-K_n$: $\forall n\geq 1, 0<\alpha_n\leq \frac{C_2}{n+C}$.

  \subsection{Estimation of $(\alpha_{-n})_{n\geq 0}$}\label{ss34}
  This time we will use the smallest fixed point $a_k=\frac{1}{2}(m_k-\sqrt{m_k^2-4})$ of $\Phi_k$ when $m_k>2$. We have (because of (\ref{E2}), $K_{-n}$ is positive):
   \begin{equation}\label{E12}   0>a_k-1>  -C_1\sqrt{m_k-2}\geq - \frac{C_2}{n+C} \end{equation}
 
  We have noticed that: $\forall n\geq 0, \beta_{-n}<K_{-n}$. Let us now compare $1+\beta_{-n-1}=\Phi_{-n}^{-1}(1+\beta_{-n})$ with $1+\beta_{-n}$. We fix a constant $B>>2$. There are three cases:
  \begin{enumerate}
  \item if $m_{-n}\leq 2$ and    $\beta_n\geq - \frac{B}{n+1+C}$; then $1+\beta_{-n-1}=\Phi_{-n}^{-1}(1+\beta_{-n})\geq \Phi_{-n}^{-1}( 1 - \frac{B}{n+1+C}))\geq  1- \frac{B}{n+1+C}$ because   $\Phi_{-n}^{-1}\geq{\rm Id}$;
   \item\label{ptbis2}  if  $m_{-n}>2$ and    $\beta_n\geq - \frac{B}{n+1+C}$;  then $1+\beta_{-n-1}=\Phi_{-n}^{-1}(1+\beta_{-n})\geq \Phi_{-n}^{-1}( 1 - \frac{B}{n+1+C}))\geq  1- \frac{B}{n+1+C}$ because $\Phi^{-1}_{-n |]-\infty, a_{-n}]}\geq {\rm Id}_{]-\infty, a_{-n}]}$ and $1-\frac{B}{n+1+C}=1- \frac{C_2}{n+1+C}-\frac{B-C_2}{n+1+C}\leq a_{-n}$ if $B$ is large enough (see (\ref{E12})); 

   \item if $\beta_{-n}< - \frac{B}{n+1+C}$; we introduce the notation $\gamma_{-n}=-\frac{B}{2(n+1+C)}$. The function $\Phi_{-n}^{-1}$ being convex such that $D(\Phi_{-n}^{-1})(1+\gamma_{-n})= \frac{1}{(m_{-n}-1-\gamma_{-n})^2}$, we have:
   $$\Phi_{-n}^{-1}(1+\gamma_{-n})-(1+\beta_{-n-1})\leq  \frac{1}{(m_{-n}-1-\gamma_{-n})^2}((1+\gamma_{-n})-(1+\beta_{-n}))$$
    If $m_{-n}\leq 2$, we have $\Phi_{-n}^{-1}(1+\gamma_{-n})\geq 1+\gamma_{-n}$ because $\Phi_{-n}^{-1}\geq {\rm Id}$; if $m_{-n}>2$, as $1+\gamma_{-n}\leq a_{-n}$ (see point \ref{ptbis2}), we have $\Phi_{-n}^{-1}(1+\gamma_{-n})\geq 1+\gamma_{-n}$ and then:
    $$\beta_{-n-1}\geq (1-\frac{1}{(m_{-n}-1-\gamma_{-n})^2})\gamma_{-n}+\frac{\beta_{-n}}{(m_{-n}-1-\gamma_{-n})^2}$$
    Because of (\ref{E8}), we have:
    \begin{equation}\label{E13}   |(m_{-n}-1-\gamma_{-n})-(1+\frac{B}{2(n+1+C)})|\leq \frac{C_2}{(n+C)^2}\end{equation} and then $ 1-\frac{1}{(m_{-n}-1-\gamma_{-n})^2}$ is positive if $C$ is large enough. Because $\gamma_{-n} >\frac{\beta_{-n}}{2}$ , we deduce:
$$\beta_{-n-1}\geq (\frac{1}{2} +\frac{1}{2(m_{-n}-1-\gamma_{-n})^2})\beta_{-n}$$and then by (\ref{E13}):
$$ \beta_{-n-1}\geq (1-\frac{B}{3(n+1+C)})\beta_{-n}
$$ If $B\geq 3$, we obtain:
\begin{equation}
\beta_{-n-1}\geq  (1-\frac{1}{(n+1+C)})\beta_{-n}\geq \frac{n+C}{n+1+C}\beta_{-n}
\end{equation}
  \end{enumerate}
   The end of the proof is then similar to the content of subsection \ref{ss33} and we obtain: $-\frac{C_2}{n+C}\leq \alpha_{-n}< 0$.

  \subsection{Regularity of the modified example}\label{ss35}
  The arguments of the proof in this subsection are very similar to the ones of M.~Herman.\\
  Let us recall that we are interested in proving that  $\displaystyle{\lim_{k\rightarrow\pm\infty} \|ÊD^2\varphi_k\|_{C^0}=0}$ and that $\displaystyle{\lim_{C\rightarrow+\infty}\sup\{  \|ÊD^2\varphi_k\|_{C^0}; k\in\Z\}=0}$
. Because of the definition $g$, we have: $\| D^2\varphi_k\|_{C^0}= \| D^2h_k+D^2h_{k-1}^{-1}-2\|_{C^0}$. Let us introduce the notation:
  $$h_k(x)=x+\Delta_k(x)=\int_0^x(1+\psi_k(t))dt.$$
  Then we want to estimate the norm $C^2$ of:
  $$\zeta_k(x)=h_k(x)+h_{k-1}^{-1}(x)-2x=\Delta_k(x)-\Delta_{k-1}(h_{k-1}^{-1}x).$$
  We differentiate to obtain:
  $$D\zeta_k(x)=D\Delta_k(x)-D\Delta_{k-1}(h_{k-1}^{-1}x)D(h_{k-1}^{-1})(x)$$
that is:
  $$D\zeta_k(x)=\psi_k(x)-\psi_{k-1}(h_{k-1}^{-1}x)D(h_{k-1}^{-1})(x).$$
  We define then $f_k:[0, \ell_k]\rightarrow [0, \ell_k]$ by: $f_k(x)=h_{k-1}(\frac{\ell_{k-1}}{\ell_k}x)$, then we have $h_{k-1}^{-1}(x)=\frac{\ell_{k-1}}{\ell_k}f_k^{-1}(x)$. We have:
  $$D(h_{k-1}^{-1})(x)=\frac{\ell_{k-1}}{\ell_k}(Df_k^{-1})(x).$$
  Let us recall that:
  $$\psi_k(x)=K_k\eta(\frac{x}{\ell_k})+\alpha_k\gamma_\pm(\frac{x}{\ell_k}).$$
Therefore
  $$\psi_{k-1}(h_{k-1}^{-1}x)=K_{k-1}\eta(\frac{h_{k-1}^{-1}x}{\ell_{k-1}})+\alpha_{k-1}\gamma_\pm(\frac{h_{k-1}^{-1}x}{\ell_{k-1}})=K_{k-1}\eta( \frac{f_k^{-1}x}{\ell_k})+\alpha_{k-1}\gamma_\pm(   \frac{f_k^{-1}x}{\ell_k}) .$$
Observe that
  $$D(h_{k-1}^{-1})(x)=\frac{1}{Dh_{k-1}(h_{k-1}^{-1}x)}=\frac{1}{1+\psi_{k-1}(h_{k-1}^{-1}x)}$$
   and then:
  $$D(h_{k-1}^{-1})(x)=\frac{1}{1+K_{k-1}\eta( \frac{f_k^{-1}x}{\ell_k})+\alpha_{k-1}\gamma_\pm(   \frac{f_k^{-1}x}{\ell_k})}
  $$
  Finally, we obtain: 
  $$\psi_{k-1}(h_{k-1}^{-1}x)D(h_{k-1}^{-1})(x)=\frac{K_{k-1}\eta( \frac{f_k^{-1}x}{\ell_k})+\alpha_{k-1}\gamma_\pm(   \frac{f_k^{-1}x}{\ell_k}) }{1+K_{k-1}\eta( \frac{f_k^{-1}x}{\ell_k})+\alpha_{k-1}\gamma_\pm(   \frac{f_k^{-1}x}{\ell_k})}.$$
 
  Moreover, we have: 
  $$Df_k(x)=\frac{\ell_{k-1}}{\ell_k}Dh_{k-1}(\frac{\ell_{k-1}}{\ell_k}x)=   \frac{\ell_{k-1}}{\ell_k}(1+   K_{k-1}\eta(\frac{x}{\ell_k})+\alpha_{k-1}\gamma_\pm(\frac{x}{\ell_k})) $$
  and
  $$Df_k^{-1}(x)= \frac{\ell_k}{\ell_{k-1}}Dh_{k-1}^{-1}(x)=\frac{\ell_k}{\ell_{k-1}}\frac{1}{1+ K_{k-1}\eta( \frac{f_k^{-1}x}{\ell_k})+\alpha_{k-1}\gamma_\pm(   \frac{f_k^{-1}x}{\ell_k})}
  $$
  Let us now compute for $x\in \T\backslash \{ \frac{\ell_k}{2}\}$ (even if $\zeta_k$ is two times differentiable at this point, the terms in the sum are not differentiable at $\frac{\ell_k}{2}$):
  $$D^2\zeta_k(x)=D\psi_k(x)-D\psi_{k-1}(h_{k-1}^{-1}x)\left( D(h_{k-1}^{-1})(x)\right)^2-\psi_{k-1}(h_{k-1}^{-1}x)D^2(h_{k-1}^{-1})(x) .$$
  Following \cite{He2}, we define:
  $$II_k= \frac{K_k}{\ell_k} D\eta(\frac{x}{\ell_k})-\frac{K_{k-1}}{\ell_k} D\eta(\frac{x}{\ell_k})+\frac{\alpha_k}{\ell_k}D\gamma_\pm(\frac{x}{\ell_k})-\frac{\alpha_{k-1}}{\ell_k}D\gamma_\pm(\frac{x}{\ell_k})$$
  $$III_k=-\left( \frac{K_{k-1}}{\ell_k} D\eta(\frac{f_k^{-1}x}{\ell_k})+\frac{\alpha_{k-1}}{\ell_k}D\gamma_\pm(\frac{f_k^{-1}x}{\ell_k})\right)\left( \frac{Df_k^{-1}(x)}{1+ K_{k-1}\eta( \frac{f_k^{-1}x}{\ell_k})+\alpha_{k-1}\gamma_\pm(   \frac{f_k^{-1}x}{\ell_k})}-1\right) $$
  $$IV_k=-\psi_{k-1}(h_{k-1}^{-1}x)D^2(h_{k-1}^{-1})(x) =\frac{\psi_{k-1}(h_{k-1}^{-1}x)D\psi_{k-1}(h_{k-1}^{-1}x)D(h_{k-1}^{-1})(x)}{(1+\psi_{k-1}(h_{k-1}^{-1}x))^2}
  $$
  i.e:
  $$IV_k=\frac{\ell_{k-1}}{\ell_k}\frac{\psi_{k-1}(h_{k-1}^{-1}x)D\psi_{k-1}(h_{k-1}^{-1}x)D(f_{k}^{-1})(x)}{(1+\psi_{k-1}(h_{k-1}^{-1}x))^2}
  $$
  and
  $$V_k=\frac{K_{k-1}}{\ell_k}\left( D\eta(\frac{x}{\ell_k})-D\eta (\frac{f_k^{-1}x}{\ell_k})\right)+
  \frac{\alpha_{k-1}}{\ell_k}\left( D\gamma_\pm (\frac{x}{\ell_k})-D\gamma_\pm (\frac{f_k^{-1}x}{\ell_k})\right).$$
  Then 
  $$D^2\zeta_k(x)=II_k+III_k+IV_k+V_k.
  $$
  Let us now estimate each term of this sum. We need some inequalities:
  \begin{equation}\label{E14}C_1\leq \| D\eta \|_{C^0},  \| D\gamma_\pm \|_{C^0},  \| D^2\eta \|_{C^0}, \| D^2\gamma_\pm \|_{C^0}\leq C_2;
  \end{equation}
 We deduce from subsections \ref{ss33} and \ref{ss34}  that:\begin{equation}\label{E15}|\alpha_k|\leq \frac{C_2}{|k|+C}\end{equation}and therefore, we have uniformly in $C>>1$ (see (\ref{E5})):
 \begin{equation}
 \lim_{k\rightarrow \pm\infty} \frac{\alpha_k^2}{\ell_k}=0.
 \end{equation} 
   From
   $$1+K_k+\alpha_k+\frac{1}{1+K_{k-1}+\alpha_{k-1}}=m_k$$
   we deduce:
   $$|\alpha_{k}-\alpha_{k-1}-(m_k-2)+K_k-K_{k-1}|\leq C_2(|K_{k-1}|+|\alpha_{k-1}| )^2$$
   and by (\ref{E3}), (\ref{E4}), (\ref{E8}), (\ref{E15}):
   \begin{equation}\label{E17} |\alpha_{k}-\alpha_{k-1}|\leq \frac{C_2}{(|k|+C)^2}
   \end{equation}
   Moreover, we have:
   $$Df_k(x)-1=  \frac{\ell_{k-1}}{\ell_k}-1+\frac{\ell_{k-1}}{\ell_k}\left(K_{k-1}\eta(\frac{x}{\ell_k})+\alpha_{k-1}\gamma_\pm(\frac{x}{\ell_k})\right)
   $$
   and 
   $$Df_k^{-1}(x)-1=  \frac{\ell_k}{\ell_{k-1}}\frac{1}{1+K_{k-1}\eta( \frac{f_k^{-1}x}{\ell_k})+\alpha_{k-1}\gamma_\pm(   \frac{f_k^{-1}x}{\ell_k})}-1
   $$
   and then we deduce from (\ref{E3}) and (\ref{E15}) that:
  \begin{equation}\label{E18}\sup\{ \| Df_k-1\|_{C^0}, \| Df_k^{-1}-1\|_{C^0}\}\leq \frac{C_2}{|k|+C}.\end{equation}
   Let us estimate $II_k=(\frac{K_k}{\ell_k} -\frac{K_{k-1}}{\ell_k} )D\eta(\frac{x}{\ell_k})+(\frac{\alpha_k}{\ell_k} -\frac{\alpha_{k-1}}{\ell_k})D\gamma_\pm(\frac{x}{\ell_k})$; because of   (\ref{E4}), (\ref{E14}) and (\ref{E17}), we have $|II_k|\leq \frac{C_2}{(|k|+C)^2\ell_k}$ and then, by (\ref{E5}),  uniformly in $C>>1$, we have:
$$ \lim_{k\rightarrow\pm \infty}|II_k|=0$$
   From (\ref{E3}), (\ref{E15}) and (\ref{E18}), we deduce that $|III_k|\leq \frac{C_2}{(|k|+C)^2\ell_k}$ and then uniformly in $C>>1$, we have:
$$\lim_{k\rightarrow\pm \infty}|III_k|=0$$  We have:
  $$IV_k=\frac{\ell_{k-1}}{\ell_k}\frac{\psi_{k-1}(h_{k-1}^{-1}x)D\psi_{k-1}(h_{k-1}^{-1}x)D(f_{k}^{-1})(x)}{(1+\psi_{k-1}(h_{k-1}^{-1}x))^2}.
  $$
  We deduce from (\ref{E3}), (\ref{E14}), (\ref{E15}) that $|IV_k|\leq \frac{C_2}{\ell_k(|k|+C)^2}$ and then that uniformly in $C>>1$, we have:
$$ \lim_{k\rightarrow\pm \infty}|IV_k|=0$$
   We have for $x\in [0, \ell_k]$: $$|D\eta(\frac{x}{\ell_k})-D\eta (\frac{f_k^{-1}x}{\ell_k})|\leq \int_0^x
   \frac{1}{\ell_k} \| D^2\eta\|_{C^0}\| Df_k^{-1}-1\|_{C^0}\leq \| D^2\eta\|_{C^0}\| Df_k^{-1}-1\|_{C^0}$$ then by  (\ref{E14}) and (\ref{E18}): $|D\eta(\frac{x}{\ell_k})-D\eta (\frac{f_k^{-1}x}{\ell_k})|\leq \frac{C_2}{|k|+C}$.  \\
Because for every $x\in [0, \ell_k]$, $\frac{x}{\ell_k}$ and $\frac{f_k^{-1}x}{\ell_k}$ are in the same half interval of $[0, 1]$, $\gamma_\pm$ is smooth between $\frac{x}{\ell_k}$ and $\frac{f_k^{-1}x}{\ell_k}$ and be can do for $\gamma_\pm$ the same estimate as for $\eta$. 
By  (\ref{E3}) and  (\ref{E15}), we deduce:

  $$|V_k|\leq \frac{C_2}{\ell_k(|k|+C)^2}$$
then, by (\ref{E5}),  uniformly in $C>>1$, we have:
$$ \lim_{k\rightarrow\pm \infty}|V_k|=0.$$   Finally, we have proved that $\varphi$ is $C^2$ and even that $\| \varphi\|_{C^2}$ is small.
  \subsection{Stable and unstable sets of the invariant curve}\label{ss36}
  We denote by $\Gamma$ the invariant curve, that is the graph of $g-{\rm Id}$.\\
  We recall that the segment with length $\frac{\ell_k}{4}$ that has same center $\mu_k$ as $I_k =[\mu_k-\frac{\ell_k}{2}, \mu_k+\frac{\ell_k}{2}]$ is denoted by $J_k=[\mu_k-\frac{\ell_k}{8}, \mu_k+\frac{\ell_k}{8}]$. Moreover, because of the definition of $h_k$ and $g_k$ (see subsection \ref{ss31}),  have: 
 \begin{enumerate}
 \item if $k\geq 1$, then: $$\forall x\in [\mu_k-\frac{\ell_k}{8}, \mu_k], g(x)=\mu_{k+1}+\frac{\ell_{k+1}}{\ell_k}(x-\mu_k);$$
 in this case, $\frac{\ell_{k+1}}{\ell_k}<1$;
 \item if $k\leq 0$, then: $$\forall x\in [\mu_k, \mu_k+  \frac{\ell_k}{8}],  g(x)=\mu_{k+1}+\frac{\ell_{k+1}}{\ell_k}(x-\mu_k);$$
if $k=0$ then  $\frac{\ell_{k+1}}{\ell_k}<1$ and if $k\leq -1$ then $\frac{\ell_{k+1}}{\ell_k}>1$.
 \end{enumerate}
 We deduce that:
 \begin{enumerate}
 \item $\forall k\geq 1, g( [\mu_k-\frac{\ell_k}{8}, \mu_k])= [\mu_{k+1}-\frac{\ell_{k+1}}{8}, \mu_{k+1}]$; then $g_{|[\mu_k-\frac{\ell_k}{8}, \mu_k]}$ is   a linear contraction;
 \item $\forall k\leq 1,   g^{-1}([\mu_k, \mu_k+  \frac{\ell_k}{8}])= [\mu_{k-1}, \mu_{k-1}+  \frac{\ell_{k-1}}{8}]$ and $(g^{-1})_{|[\mu_k, \mu_k+  \frac{\ell_k}{8}]}$ is linear, a contraction if $k\leq 0$ and a dilatation if $k=1$.
 \end{enumerate}
 We introduce the family $(S_k)_{k\geq 1}$ and $(U_k)_{k\leq 0}$ of segments of $\T\times \R$ defined by:
 $$S_k=\{ (x, g(x)-x); x\in  [\mu_k-\frac{\ell_k}{8}, \mu_k]\}\quad{\rm and}\quad U_k=\{ (x, g(x)-x); x\in [\mu_k, \mu_k+  \frac{\ell_k}{8}]\}.$$
 Because the curve $\Gamma$ is the graph of $g-{\rm Id}$, these segments are subsets of $\Gamma$.
 We have: 
 $$\forall k\geq 1, f_\varphi(S_k)=S_{k+1}\quad{\rm and}\quad \forall k\leq 0, f_\varphi^{-1}(U_k)=U_{k-1};$$
 in the first case, $f_{\varphi|S_k}$ is a linear contraction with rapport $\frac{\ell_{k+1}}{\ell_k}$ and in the second case  $f^{-1}_{\varphi|U_k}$ is a linear contraction with rapport $\frac{\ell_{k-1}}{\ell_{k}}$. \\

We have proved in subsection \ref{ss32} some equalities for $h_k+h_{k-1}^{-1}$ that implies:
$$\forall k\in \Z, \forall x\in J_k, \varphi (x)= (m_k-2)(x-\mu_k)+\mu_{k+1}+\mu_{k-1}-2\mu_k.$$
Let us recall that:
$$f_\varphi(\theta,r)=(\theta+r, r+\varphi(\theta+r))\quad{\rm and}\quad f_\varphi^{-1}(\theta, r)=(\theta-r+\varphi (\theta), r-\varphi(\theta))
$$
therefore the restriction of  $f_\varphi^{-1}$ to any band $J_k\times \R$ is linear. If we know the expression of a linear map on a segment, we can deduce the expression of the map on the whole line supporting the segment. In particular, if we define the families of segments $(\tilde S_k)_{k\geq 1}$ and $(\tilde U_k)_{k\leq 0}$ by :
\begin{equation}\label{E20}\tilde S_k=\{ (x, \mu_{k+1}-\mu_k+(\frac{\ell_{k+1}}{\ell_k}-1)(x-\mu_k)); x\in J_k\}\quad{\rm for}\quad k\geq 1\end{equation} and \begin{equation} \tilde U_k=\{ (x, \mu_{k+1}-\mu_k+(\frac{\ell_{k+1}}{\ell_k}-1)(x-\mu_k)); x\in J_k\}\quad{\rm for}\quad k\leq 0\end{equation}
then we have $U_k\subset \tilde U_k$, $S_k\subset \tilde S_k$ and:
$$\forall k\geq 1, f_\varphi(\tilde S_k)=\tilde S_{k+1}\quad{\rm and}\quad \forall k\leq 0, f_\varphi^{-1}(\tilde U_{k})=\tilde U_{k-1}.$$
 Moreover, the restriction of $f_\varphi$ to $\tilde S_k$ is a linear contraction with rapport $\frac{\ell_{k+1}}{\ell_k}$ and the restriction of $f_\varphi^{-1}$ to $\tilde U_k$ is a linear contraction with rapport $\frac{\ell_{k-1}}{\ell_{k}}$. We then deduce that $\tilde S_k$ is in the stable set of the point $(\mu_k, \mu_{k+1}-\mu_k)$ and that $\tilde U_k$ is in the unstable set of the point $(\mu_k, \mu_{k+1}-\mu_k)$.
 
 We then extend these two families of segments by:
 \begin{enumerate}
 \item[$\bullet$] if $k\leq 0$, $\tilde S_k=f_\varphi^{k-1}(\tilde S_1)$;
 \item[$\bullet$] if $k\geq 0$, $\tilde U_k=f_\varphi^{k}(\tilde U_0)$.
 \end{enumerate}

  Let us now choose a $C^\infty$ injective map $\gamma_1^s: \R\rightarrow \tilde S_1$ such that $\gamma_1^s(0)=x_1=(\mu_1, \mu_2-\mu_1)$, $\gamma_1^s(\R)$ is $\tilde S_1$ without its ends and  $\gamma_1^s(]-\infty, 0[)$ is $S_1$ without its ends. 
 \\
 Similarly, we choose a $C^\infty$ injective map $\gamma_0^u: \R\rightarrow \tilde U_0$ such that $\gamma_0^u(0)=x_0=(\mu_0, \mu_1-\mu_0)$,   $\gamma_0^u(\R)$ is $\tilde U_0$ without its ends and  $\gamma_0^u(]0,+\infty[)$ is $U_0$ without its ends. 
 
 We   extend these curves to two families by: $\gamma_k^s=f_\varphi^{k-1}\circ \gamma_1^s$ and $\gamma_k^u=f_\varphi^k\circ \gamma_0^u$.  Then we have: 
\begin{enumerate}
\item $f_\varphi\circ \gamma_k^u=\gamma_{k+1}^u$ and $f_\varphi\circ \gamma_k^s=\gamma_{k+1}^s$;
\item $\forall y\in \gamma_0^s (\R) , \displaystyle{\lim_{n\rightarrow +\infty}d(f_\varphi^ny, f_\varphi^nx_0)=0}$ and 
$\forall y\in  \gamma_0^u (\R), \displaystyle{\lim_{n\rightarrow +\infty}d(f_\varphi^{-n}y, f_\varphi^{-n}x_0)=0}$;
\item  $\gamma_k^s(]-\infty, 0])\cup \gamma_k^u([0, +\infty[)\subset f_\varphi^{k-1}(S_1)\cup f_\varphi^k(U_0)\subset \Gamma$.
\end{enumerate}

Let us now prove that $\gamma_1^s(]0, +\infty[)\cup \gamma_0^u(]-\infty, 0[)\subset \T\times \R\backslash \Gamma$. We will deduce that $\gamma_1^s(]0, +\infty[)$ is a part of the stable set of $x_1$ and  of $\Gamma$ that doesn't meet $\Gamma$, hence it is in an instability zone $\Uc$ and $\Gamma$ is in the boundary of $\Uc$ (see the proposition contained in subsection \ref{ss13}); we will even see that $\Uc$ is under $\Gamma$. Similarly, we will prove that $\gamma_0^u(]-\infty, 0[)$ is in $\Uc$. We will of course deduce that:
$$\forall k\in \Z, \gamma_k^s(]-\infty, 0[)\cup\gamma_k^u(]0, +\infty[)\subset \Uc.$$

By (\ref{E20}), we  have an explicit expression for:
$$\gamma_1^s(]0, +\infty[)=\{(x, \mu_2-\mu_1+(\frac{\ell_2}{\ell_1}-1)(x-\mu_1)); x\in ]\mu_1, \mu_1+\frac{\ell_1}{8}[\}$$
Moreover, because of the definition of $\Gamma$, we have (see subsection \ref{ss31}):
$$\forall x\in  ]\mu_1, \mu_1+\frac{\ell_1}{8}[, g(x)-x=(\frac{\ell_2}{\ell_1}+\alpha_1-1)(x-\mu_1)+\mu_2-\mu_1.$$
As $\alpha_1>0$ and $\Gamma$ is the graph of $g-{\rm Id}$, we deduce that $\gamma_1^s(]0, +\infty[)$ doesn't meet $\Gamma$, and even that $\gamma_1^s(]0, +\infty[)$  is under $\Gamma$. A similar argument gives the result for $\gamma_0^u(]-\infty, 0[)$.

\begin{remk}
1) If we exchange $\gamma_-$ and $\gamma_+$, we obtain an instability zone $\Uc$ that is above $\Gamma$.\\
2) If we use a similar construction along two wandering intervals, we obtain a curve $\Gamma$ that is at the boundary of two instability zones. 
\end{remk}

\section{The case of the $C^1$ topology: proof of theorem \ref{T2}\!}\label{S4}
If $U$ is an open subset of $\A$, the set ${\rm Diff}^1_\omega (U)$ of $C^1$  symplectic diffeomorphisms of $U$ is endowed with the  strong Whitney's topology  (see \cite{Hir, PR}).  Observe that the set $\Tc$ of symplectic twist maps of $\A$ is open for the Whitney topology in  ${\rm Diff}^1_\omega (\A)$. Then if $U$ is any open subset of $\A$, the set $\Tc(U)$ of the restrictions to $U$ of symplectic twist maps of $\A$ is open in  ${\rm Diff}^1_\omega (U)$.\\
The following result is theorem 3 of  \cite{ABC}.
\begin{theo} {\bf (\cite{ABC})} Let $(M, \omega)$ be a non-compact closed manifold. There exists a dense $G_\delta$ subset $\Gc$ of ${\rm Diff}^1_\omega (M)$ such that, for all $f\in\Gc$, the set of points of $M$ whose positive orbit is relatively compact in $M$ has no interior.

\end{theo}
Let us now consider an essential invariant curve $\Gamma$ of a symplectic twist map $f$ of $\A$. The curve $\Gamma$ is then the graph of a Lipschitz map $\gamma: \T\rightarrow \R$. Denoting by $U$ one of the two connected components of $\A\backslash \Gamma$, we have: $f(U)=U$. In order to define a neighborhood $\Uc$ of $f_{|U}$ for the $C^1$ strong topology, we use the function $\varepsilon: U\rightarrow \R_+^*$ defined by $\varepsilon (\theta, r)=(r-\gamma(\theta))^2$.
$$\Uc=\{ g\in \Tc (U); \forall (\theta, r)\in U, \sup\left\{ d(f(\theta, r), g(\theta, r)),  \| Df(\theta, r)-Dg(\theta, r)\|\right\} \leq \varepsilon (\theta, r)\}.$$
The previous theorem  implies that there exists $h\in \Uc$ such that 
   the set of points of $U$ whose positive orbit for $h$ is relatively compact in $U$ has no interior.

We now define $g: \A\rightarrow \A$ such that $g_{|\A\backslash U}=f_{|\A\backslash U}$ and $g_{|U}=h$. It comes from the definition of $\Uc$ and  $\Tc (U)$  that $g\in {\rm Diff}^1_\omega (\A)$.\\
Let us prove that $U$  contains at most one essential invariant curve for $g$. If not, there exists  a bounded invariant open region $R$ between two such invariant curves. Then all the points of $R$ have a positive orbit that is relatively compact in $U$, this is a contradiction with the choice of $h$. \\
Let us now assume that $f_{|\Gamma}$ has an irrational rotation number or that $f_{|\Gamma}$ is $C^0$ conjugated to a rational rotation. We have noticed in subsection \ref{ss13} that in this case, if $\Gamma^*$ is another essential invariant curve of $f$, than $\Gamma\cap \Gamma^*=\emptyset$. Hence the closure $\bar U$ of $U$ contains at most one essential invariant curve that is different from $\Gamma$, and this curve is contained in $U$. There are two cases: 
\begin{enumerate}
\item[$\bullet$] either $U$ contains one essential invariant curve $\Gamma'$ for $g$.  The region $R$ between $\Gamma$ and $\Gamma'$ is an instability zone for $g$ and its boundary contains $\Gamma$;
\item[$\bullet$] or $U$ contains no essential invariant curve for $g$. The curve $\Gamma$ is at the boundary of the instability zone $U$ of $g$.
\end{enumerate}

\begin{remk}
Using   methods contained in the (non-published) thesis of my student Marie Girard that allows us to destroy all the invariant curves by perturbation, we can choose $h$ such that $U$ contains no essential invariant curve and thus is an instability zone. 
\end{remk}
The only case that we did not solve is the case of a rational rotation number for $f_{|\Gamma}$  and a hyperbolic dynamics for $f_{|\Gamma}$ (i.e. we assume that all the periodic points of $f_{|\Gamma}$ are hyperbolic).\\
 C.~Robinson proved the following result in \cite{Rob}. Let us recall that a periodic point $p$ of $f$ with period $\tau$ is non-degenerate if no root of  1 is an eigenvalue of $Df^\tau (p)$.
\begin{theo}{\bf(\cite{Rob})} Let $(M, \omega)$ be a  closed manifold. There exists a dense $G_\delta$ subset $\Gc$ of ${\rm Diff}^1_\omega (M)$ such that, for all $f\in\Gc$, the periodic points are non degenerate and the stable and unstable manifolds of each pair of hyperbolic periodic orbits of $f$ are transverse at all of their points  of intersection. 

\end{theo}

We assume  now that $f\in{\Tc }(A)$ a symplectic twist map that has an essential invariant curve $\Gamma$ such that:
\begin{enumerate}
\item[$\bullet$] its rotation number is rational;
\item[$\bullet$] the periodic points of $f_{|\Gamma}$ are all hyperbolic.
\end{enumerate}
Then there is a finite number of such periodic points, that we denote by $x_1, \dots, x_n$, and $\Gamma$ is the union of $\{ x_1, \dots , x_n\}$ and some branches of the stable/unstable manifolds of these periodic points. Let us notice that every $g\in\Uc$ can be extend in a unique $\tilde g\in \Tc(\A)$ by: $\tilde g_{|\A\backslash \Uc}=f_{|\A\backslash \Uc}$ and that in this case, $D\tilde g_{|\Gamma}=Df_{|\Gamma}$. Hence  $\tilde g$ has the same periodic points as $f$ on $\Gamma$, and this periodic points are hyperbolic. 

 We can directly adapt Robinson's proof to build a dense $G_\delta$ $\Gc$ of $\Uc$ such that for all $g\in \Gc$, the stable and unstable branches of the stable  and unstable  manifolds of the $x_i$ for $\tilde g$ that are contained in $U$ are transverse at all of their points of intersection. 
 
 If now $\Gamma^*$ is an essential invariant curve for $\tilde g$ that is contained in $\bar U$ and that meets $\Gamma$, then $\Gamma\cap \Gamma^*$ is a closed invariant set that contains a point of the stable manifold of a point  $x_i$. Hence it contains this $x_i$. The rotation number of $\Gamma^*$ is then equal to the one of $\Gamma$, and then $\Gamma^*$ is  the union of $\{ x_1, \dots , x_n\}$ and some branches of the stable/unstable manifolds of these periodic points.  But if $\Gamma\not=\Gamma^*$, then $\Gamma^*$ contains a branch in $U$ that is a stable and an unstable branch, and this contradicts the transversality  of such branches. We deduce that either $\Gamma=\Gamma^*$ or $\Gamma\cap \Gamma^*=\emptyset$, and we can conclude exactly in the same way as is the irrational case.

\end{document}